\documentclass[a4paper]{article}
\usepackage[utf8]{inputenc}
\usepackage{latexsym} 
\usepackage{amsmath,color}
\usepackage{amsfonts}
\usepackage{amssymb}
\usepackage{mathtools}
\usepackage{cite}
\newlength{\bibitemsep}
\newlength{\bibparskip}
\let\oldthebibliography\thebibliography
\renewcommand\thebibliography[1]{%
  \oldthebibliography{#1}%
  \setlength{\parskip}{\bibitemsep}%
  \setlength{\itemsep}{\bibparskip}%
}
\usepackage{amsthm}
\newtheorem{thm}{Theorem}%[within]
\newtheorem{cor}{Corollary}%[within]
\newtheorem{rem}{Remark}%[within]

\renewcommand{\Re}{\operatorname{Re}}

\usepackage{graphicx}
\usepackage{multicol}
\newcommand{\msc}[1]{\ \footnote{\hspace{-7mm}
\rm 2010 Mathematics Subject Classification: #1}}
\newcommand{\keywords}[1]{\ \footnote{\hspace{-7mm}
\rm Key Words and Phrases: \it #1}}

\makeatletter

\title{The limit of the Riemann zeta function and its nontrivial zeros}%%%%[]

\author{Tanfer Tanriverdi\\
Department of Mathematics  Faculty of Arts and Sciences\\ 
Harran University \c{S}anl{\i}urfa Turkey 63290\\
\\
 ttanriverdi@harran.edu.tr
\\
\\
}

\date{}
\begin{document}
\maketitle
\begin{abstract}
In this article, with a new approach, which is not discussed in the literature yet, the limit of the Riemann zeta function or Euler-Riemann zeta function is approximately explored by applying Dirichlet's rearrangement theorem for absolutely convergent series to the Riemann zeta function by rearranging its terms as geometric series for sufficiently large $n$. The limit of the Riemann zeta function or Euler-Riemann zeta functions, $\lim_{n\to\infty} \zeta(z)$, is first time explored.  The limit obtained here is a very promising for the nontrivial or complex zeros of  the Rieman zeta function.
\end{abstract}
\msc{11M06, 11M26, 40A25, 30B50} 
%%%%%%%%%%%%%%%%%%%%%%%%%%%%%%%%
\section{Introduction }\label{intro}
%%% add text here for the Introduction
The importance of the Riemann zeta function is trivial to anybody who runs into it while working out with certain subjects such as number theory for investigating properties of prime numbers, complex analysis, physics etc. Riemann zeta function and its crucial properties have been studied extensively, see \cite{1,2,3,4,5,6,7,8,9,10,11,12}.

The closed form formula for the values $\zeta(2m)$ in terms of Bernoulli number is well known due to  Euler  where $m$ is an integral \cite{13}. The values $\zeta(2m+1)$ are not known yet except $\zeta(3)$  Ap{\'e}ry \cite{14} however some representations of $\zeta(2m+1)$ are given in terms of special functions. Several interesting  evaluations and representations of the Riemann zeta function $\zeta(m)$ with $m\ne 1$ are reviewed in \cite{9}.

It was conjectured by Riemann \cite{7} that all nontrivial zeros of  $\zeta(z)$ lie on the line $\Re(z)=\frac{1}{2}$ but this hypothesis is never proved or disapproved. Theoric and numerical study were carried out strongly by \cite{1,2} and references therein.

We also aim to say something about the nontrivial zeros of the Riemann zeta function. To reach this goal we need to find the limit of the Riemann zeta function to see how zeros of the its limit behaves. We here apply Dirichlet's rearrangement theorem for absolutely convergent series to the Riemann zeta function by rearranging its terms. Such action is legal since the Riemann zeta function is also convergent absolutely where $\Re(z)>1$. Then, one rewrites Riemann zeta function by rearranging its terms as geometric series as given below  for sufficiently large $n$ to get its limit involving infinite series with the hyperbolic functions have been also attracted the attention of many authors \cite{15,16}.

Therefore, we are about to present a very interesting since it is new and elementary approach to Riemann zeta function to get its limit that may attract to interested readers. This is a new way of looking at Riemann zeta function in terms of calculating its limit.  

In Section  \ref{some},  some  classical variations of the Riemann zeta function and in Section \ref{mainresult}, main results are introduced. In final Section \ref{con}, results ontained in this paper are highlighted.
%%%%%%%%%%%%%%%%%%%%%%%%%%%%
\section{Some variations of the Riemann zeta function}  \label{some}   
In Riemann \cite{7}, the Riemann zeta function $\zeta(z)$ is defined by 
\begin{eqnarray}\label{eq1}
\zeta(z)=\sum_{n=1}^{\infty} \frac{1}{n^{z}}, \ \Re(z)>1.
\end{eqnarray}
The series (\ref{eq1}) is absolutely convergent and analytic function for $\Re(z)>1$. 
If $z=s>1$ is real number, 
then Riemann zeta function is represented and studied by Euler as 
\begin{eqnarray*}
\zeta(z)=\prod_{p} \left(1-\frac{1}{p^{s}}\right)^{-1}
\end{eqnarray*}
where $p$ runs through all prime numbers. There are several representations for Rieman zeta function  
\begin{eqnarray*}
\zeta(z)= \frac{1}{\Gamma(z)}\int_{0}^{\infty}\frac{x^{z-1}}{e^{x}-1}dx, \ \Re(z)>1
\end{eqnarray*}
and 
\begin{eqnarray*}\label{eq2}
\zeta(z)= 2^{z}\pi^{z-1}\sin(\frac{z\pi}{2})\Gamma(1-z)\zeta(1-z), \ \Re(z)>1.
\end{eqnarray*}
One may extend analytic region of Riemann zeta function from $\Re(z)>1$ to $\Re(z)>0$.
\begin{eqnarray}\label{eq3}
\left(1-2^{1-z}\right)\zeta(z)=\sum_{n=1}^{\infty}\frac{(-1)^{n-1}}{n^{z}}, \ \Re(z)>0 \ \mbox{and}\ z\neq 1.
\end{eqnarray}
For more classical representations of the Riemann zeta function, see \cite{1,2}. 
%%%%%%%%%%%%%%%%%%%%%%%%%%%%%%%%%%%%%%%%%%%%%%%%%%%%%%%%%%%%%%%%%
\section{Main results} \label{mainresult}
Since Riemann zeta function is convergent absolutely when $\Re(z)>1$. Then, one rewrites Riemann zeta function by rearranging its terms as the following for sufficiently large $n$.
\begin{thm} \label{thm1}
Let (\ref{eq1}) hold. Then limit of the Riemann zeta function is given by for sufficiently large $n$ as
\begin{eqnarray*}\label{eq4}
\begin{aligned}
\zeta(z)&=1+\frac{1}{2}\overbrace{\sum_{r=2}^{n}\coth\left(\frac{z\log r}{2}\right)}^{\text{$l=$ number of terms}}-\frac{l}{2}\\
&=\frac{2-l}{2}+\frac{1}{2}\sum_{r=2}^{n}\coth\left(\frac{z\log r}{2}\right)\\
\end{aligned}
\end{eqnarray*}
where $r\neq 2^{k},3^{k},5^{k},6^{k},7^{k},10^{k},11^{k},12^{k},\cdots $ with $k=2,3,4,\cdots$. 
\end{thm}
\begin{proof}
For sufficiently large $n$, one rewrites Riemann zeta function (\ref{eq1}) by rearranging its terms as
\begin{eqnarray}\label{eq5}
\begin{aligned}
\zeta(z)&=1+\frac{1}{2^{z}}+\frac{1}{2^{2z}}+\frac{1}{2^{3z}}+\frac{1}{2^{4z}}+\cdots\\
&+\frac{1}{3^{z}}+\frac{1}{3^{2z}}+\frac{1}{3^{3z}}+\frac{1}{3^{4z}}+\cdots \\
&+\frac{1}{5^{z}}+\frac{1}{5^{2z}}+\frac{1}{5^{3z}}+\frac{1}{5^{4z}}+\cdots\\
&+\frac{1}{6^{z}}+\frac{1}{6^{2z}}+\frac{1}{6^{3z}}+\frac{1}{6^{4z}}+\cdots\\
%%&+\frac{1}{7^{z}}+\frac{1}{7^{2z}}+\frac{1}{7^{3z}}+\frac{1}{7^{4z}}+\cdots\\
%%&+\frac{1}{10^{z}}+\frac{1}{10^{2z}}+\frac{1}{10^{3z}}+\frac{1}{10^{4z}}+\cdots\\
&+\cdots+\frac{1}{n^{z}}+\frac{1}{n^{2z}}+\frac{1}{n^{3z}}+\frac{1}{n^{4z}}+R_{n}(z) 
\end{aligned}
\end{eqnarray}
where $n\neq 2^{k},3^{k},5^{k},6^{k},7^{k},10^{k},11^{k},12^{k},\cdots $ with  $k=2,3,4,\cdots$ and $R_{n}(z)$ is the remainder. It is important to note that these excluded values are already included in (\ref{eq5}). So by applying geometric series to (\ref{eq5}) for sufficiently large $n$, one obtains 
the approximate limit to the Riemann zeta function as
\begin{eqnarray}\label{eq6}
\begin{aligned}
\zeta(z)&= \frac{1}{1-\frac{1}{2^{z}}}+\left(\frac{1}{1-\frac{1}{3^{z}}}-1\right)+\left(\frac{1}{1-\frac{1}{5^{z}}}-1\right) 
+\left(\frac{1}{1-\frac{1}{6^{z}}}-1\right)\\
&+\left(\frac{1}{1-\frac{1}{7^{z}}}-1\right)+\left(\frac{1}{1-\frac{1}{10^{z}}}-1\right)+\left(\frac{1}{1-\frac{1}{11^{z}}}-1\right)\\
&+\left(\frac{1}{1-\frac{1}{12^{z}}}-1\right)+\cdots+\left(\frac{1}{1-\frac{1}{n^{z}}}-1\right).\\
\end{aligned}
\end{eqnarray}
Here, $R_{n}(z)$ goes to zero for sufficiently large $n$. Indeed, if $\Re(z)=\sigma>1$ then
\begin{eqnarray}\label{eq7}\nonumber
\lvert R_{n}(z)\rvert\leq \frac{1}{(n+1)^{\sigma}}+ \frac{1}{(n+2)^{\sigma}}+\frac{1}{(n+3)^{\sigma}}+\cdots. 
\end{eqnarray}
So, $ \lvert R_{n}(z)\rvert \to 0$ for $n \to \infty$.
After arranging (\ref{eq6}), one obtains the approximate limit as
\begin{eqnarray}\label{eq8}
\begin{aligned}
\zeta(z)&= 1+\frac{1}{2^{z}-1}+\frac{1}{3^{z}-1}+\frac{1}{5^{z}-1}+\frac{1}{6^{z}-1}+\frac{1}{7^{z}-1}\\
&+\frac{1}{10^{z}-1}+\frac{1}{11^{z}-1}+\frac{1}{12^{z}-1}+\cdots+\frac{1}{n^{z}-1}.\\
\end{aligned}
\end{eqnarray}
One may rewrite (\ref{eq8}) as series  
\begin{eqnarray}\label{eq9}
\begin{aligned}
\zeta(z)=1+\sum_{r=2}^{n}\frac{1}{r^{z}-1}
\end{aligned}
\end{eqnarray}
where $r\neq 2^{k},3^{k},5^{k},6^{k},7^{k},10^{k},11^{k},12^{k},\cdots $ with $k=2,3,4,\cdots$.  
\begin{eqnarray}\label{eq10}
\frac{1}{e^{z}-1}+\frac{1}{2}=\frac{1}{2}\coth\left(\frac{z}{2}\right).
\end{eqnarray}
It is well known that equation (\ref{eq10}) has a very important role in Riemann zeta function.
So one also rewrites (\ref{eq9}) by using  (\ref{eq10}) as 
\begin{eqnarray}\label{eq11}
\begin{aligned}
\zeta(z)&=1+\frac{1}{2}\overbrace{\sum_{r=2}^{n}\coth\left(\frac{z\log r}{2}\right)}^{\text{$l=$ number of terms}}-\frac{l}{2}\\
&=\frac{2-l}{2}+\frac{1}{2}\sum_{r=2}^{n}\coth\left(\frac{z\log r}{2}\right)\\
\end{aligned}
\end{eqnarray}
where $r\neq 2^{k},3^{k},5^{k},6^{k},7^{k},10^{k},11^{k},12^{k},\cdots $ with $k=2,3,4,\cdots$. 
\end{proof}
Infinite series involving the hyperbolic functions have attracted the attention of many authors \cite{15,16}.
%%%%%%%%%%%%%%%%%%%%%%%%%%%%%%%%%%%%%%%%%%%%%%%%%%%%%%%%%%%%%%%%
\begin{rem} \label{rem1}
Let $r\neq 2^{k},3^{k},5^{k},6^{k},7^{k},10^{k},11^{k},12^{k},\cdots $ be where $k=2,3,4,\cdots$. Then, it is easy to see that (\ref{eq5}) has the following series representation.
\begin{eqnarray}\label{eq8a}
\begin{aligned}
\zeta(z)&=\sum_{n=-1}^{\infty}\frac{z^{n}B_{1+n}\sum_{r=2}^{n}\log^{n} r}{(1+n)!}\\
&=\sum_{n=-1}^{\infty}\frac{z^{n}B_{1+n}\left(\log^{n} 2+\log^{n} 3+\cdots+\log^{n} n\right)}{(1+n)!}.
\end{aligned}
\end{eqnarray}
\end{rem} 
Here, $ B_{1+n}$ is $(1+n)$ th Bernoulli numbers \cite{12} and $r$ is defined as above. 
%%%%%%%%%%%%%%%%%%%%%%%%%%%%%%%%%%%%
From (\ref{eq8}) and (\ref{eq11}), one gets the following limiting approximations for zeta function where $m$ is a positive integer.
\begin{cor}\label{cor1}
Let $r\neq 2^{k},3^{k},5^{k},6^{k},7^{k},10^{k},11^{k},12^{k},\cdots $ be where $k=2,3,4,\cdots$. Then
\begin{eqnarray*}\label{eq12}
\begin{aligned}
\zeta(m)&=1+\overbrace{\sum_{r=2}^{n} \frac{1}{r^{m}-1}}^{\text{$l=$ number of terms}}\\
&=\frac{2-l}{2}+\frac{1}{2}\sum_{r=2}^{n}\coth\left(\frac{m\log r)}{2}\right).
\end{aligned}
\end{eqnarray*}
\end{cor}
\begin{proof}
Proof is the same as in Theorem \ref{thm1} when $z$ is replaced by $m$.
\end{proof}
\begin{cor} \label{cor2}
Let $r\neq 2^{k},3^{k},5^{k},6^{k},7^{k},10^{k},11^{k},12^{k},\cdots $ be where $k=2,3,4,\cdots$. Then
\begin{eqnarray*}\label{eq13}
\begin{aligned}
\zeta(2m)&=1+\overbrace{\sum_{r=2}^{n} \frac{1}{r^{2m}-1}}^{\text{$l=$ number of terms}}\\
&=\frac{2-l}{2}+\frac{1}{2}\sum_{r=2}^{n}\coth\left(m\log r\right).
\end{aligned}
\end{eqnarray*}
\end{cor}
\begin{proof}
Proof is the same as in Theorem \ref{thm1} when $z$ is replaced by $2m$.
\end{proof}
The values $\zeta(2m)$ are well known due to Euler where $m$ is an integral \cite{13}.
\begin{cor} \label{cor3}
Let $r\neq 2^{k},3^{k},5^{k},6^{k},7^{k},10^{k},11^{k},12^{k},\cdots $ be where $k=2,3,4,\cdots$. Then
\begin{eqnarray*}\label{eq14}
\begin{aligned}
\zeta(2m+1)&=1+\overbrace{\sum_{r=2}^{n} \frac{1}{r^{2m+1}-1}}^{\text{$l=$ number of terms}}\\
&=\frac{2-l}{2}+\frac{1}{2}\sum_{r=2}^{n}\coth\left((m+\frac{1}{2})\log r\right).
\end{aligned}
\end{eqnarray*}
\end{cor}
\begin{proof}
Proof is the same as in Theorem \ref{thm1} when $z$ is replaced by $2m+1$.
\end{proof}

The values $\zeta(2m+1)$ are not known yet except $\zeta(3)$ \cite{14} however some representations of $\zeta(2m+1)$ are given in terms of special functions. 

%%One may want to know the zeros of the  approximate limits obtained above.  We may at least reach this goal by using the program in \cite{14} for some %%terms. By using "solve" command, one obtains the followings.
%%\begin{eqnarray*}
%%1+\frac{1}{2^x-1}=0
%%\end{eqnarray*}
%%There are not any solutions. It is also easy to see that without using any program since it is monotonically decreasing. 
%%\begin{eqnarray*}
%%1+\frac{1}{2^x-1}+\frac{1}{3^x-1}=0
%%\end{eqnarray*}
%%There are two complex solutions $ x\approx -5.30989\times10^{-17}\pm 3.50671i\dots{}$   
%%\begin{eqnarray*}
%%1+\frac{1}{2^x-1}+\frac{1}{3^x-1}+\frac{1}{5^x-1}=0
%%\end{eqnarray*}
%%There are two complex solutions $ x\approx 0.445959\pm 2.81436i\dots{}$ 
%%\begin{eqnarray*}
%%1+\frac{1}{2^x-1}+\frac{1}{3^x-1}+\frac{1}{5^x-1}+\frac{1}{6^x-1}=0
%%\end{eqnarray*}
%%There are two complex solutions $ x\approx 0.631214\pm 2.54663i\dots{}$.  

By applying the same argument as in (\ref{eq5}) to  (\ref{eq3}), Then one formally obtains the following theorem.
\begin{thm} \label{thm2} 
Let (\ref{eq3}) hold. Then the approximate limit of the Riemann zeta function is given by for sufficiently large $n$ as
\begin{eqnarray*}\label{eq15}
\zeta(z)=\begin{cases}\frac{1}{1-2^{1-z}}\left(1+\frac{1}{2}\overbrace{\sum_{r=2}^{n}(-1)^{r-1}\coth\left(\frac{z\log r}{2}\right)}^{\text{$l=$ number of terms is even}}\right), & \text{if $n$ is even}\\
\frac{1}{1-2^{1-z}}\left(\frac{1}{2}+\frac{1}{2}\overbrace{\sum_{r=2}^{n}(-1)^{r-1}\coth\left(\frac{z\log r}{2}\right)}^{\text{$l=$ number of terms is odd}}\right), 
& \text{if $n$ is odd}.\\
\end{cases}
\end{eqnarray*}
\end{thm}
\begin{rem} Similar Corollaries  may be derived related to $\zeta(m)$,  $\zeta(2m)$ and  $\zeta(2m+1)$ as stated in Corollary \ref{cor1}, \ref{cor2} and \ref{cor3}.
\end{rem}
\begin{proof}
One rewrites (\ref{eq3}) for sufficiently large $n$ by rearranging its terms as
\begin{eqnarray}\label{eq16}
\begin{aligned}
\left(1-2^{1-z}\right)\zeta(z)&=\Big(1+\frac{1}{3^{z}}+\frac{1}{3^{2z}}+\frac{1}{3^{3z}}+\frac{1}{3^{4z}}+\cdots\\
&+\frac{1}{5^{z}}+\frac{1}{5^{2z}}+\frac{1}{5^{3z}}+\frac{1}{5^{4z}}+\cdots \\
&+\frac{1}{7^{z}}+\frac{1}{7^{2z}}+\frac{1}{7^{3z}}+\frac{1}{7^{4z}}+\cdots\Big)\\
%%%& +\frac{1}{11^{z}}+\frac{1}{11^{2z}}+\frac{1}{11^{3z}}+\frac{1}{13^{4z}}+\cdots\\
%%%&+\frac{1}{13^{z}}+\frac{1}{13^{2z}}+\frac{1}{13^{3z}}+\frac{1}{13^{4z}}+\cdots\Big)\\
&+\cdots-\Big(\frac{1}{2^{z}}+\frac{1}{2^{2z}}+\frac{1}{2^{3z}}+\frac{1}{2^{4z}}+\cdots\\
&+\frac{1}{6^{z}}+\frac{1}{6^{2z}}+\frac{1}{6^{3z}}+\frac{1}{6^{4z}}+\cdots\\
&+\frac{1}{10^{z}}+\frac{1}{10^{2z}}+\frac{1}{10^{3z}}+\frac{1}{10^{4z}}+\cdots\Big)\\
%%%&+\frac{1}{12^{z}}+\frac{1}{12^{2z}}+\frac{1}{12^{3z}}+\frac{1}{12^{4z}}+\cdots\\
%%%&+\frac{1}{14^{z}}+\frac{1}{14^{2z}}+\frac{1}{14^{3z}}+\frac{1}{14^{4z}}+\cdots\Big)+\cdots\\ 
%&+\cdots+\frac{1}{n^{z}}+\frac{1}{n^{2z}}+\frac{1}{n^{3z}}+\frac{1}{n^{4z}}+R_{n}(z)
\end{aligned}
\end{eqnarray}
%%%%where $n\neq 2^{k},3^{k},5^{k},6^{k},7^{k},10^{k},11^{k},12^{k},\cdots$ and  $k=2,3,4,\cdots$. 
So by applying geometric series to (\ref{eq16}) as before one obtains 
\begin{eqnarray}\label{eq17}
\begin{aligned}
\left(1-2^{1-z}\right)\zeta(z)&= 1+\frac{1}{3^{z}-1}-\frac{1}{2^{z}-1}+\frac{1}{5^{z}-1}-\frac{1}{6^{z}-1}\\
&+\frac{1}{7^{z}-1}-\frac{1}{10^{z}-1}+\frac{1}{11^{z}-1}-\frac{1}{12^{z}-1}\\
&+\frac{1}{13^{z}-1}-\frac{1}{14^{z}-1}+\cdots +R_{n}(z).
%%%%+\sum_{r}^{n}(-1)^{r-1}\frac{1}{r^{z}-1}\\
\end{aligned}
\end{eqnarray}
As before, $ \lvert R_{n}(z)\rvert \to 0$ as $n \to \infty$. One may rewrite   approximation to  (\ref{eq17}) as
\begin{eqnarray}\label{eq18}
\begin{aligned}
\zeta(z)=\frac{1}{1-2^{1-z}}\left(1+\sum_{r=2}^{n}(-1)^{r-1}\frac{1}{r^{z}-1}\right)
\end{aligned}
\end{eqnarray}
where $r\neq 2^{k},3^{k},5^{k},6^{k},7^{k},10^{k},11^{k},12^{k},\cdots $ and $k=2,3,4,\cdots$. 
Once again by using  (\ref{eq10}), one obtains the following limiting representation to Riemann zeta function for sufficiently large $n$ as
\begin{eqnarray*}\label{eq19}
\zeta(z)=\begin{cases}\frac{1}{1-2^{1-z}}\left(1+\frac{1}{2}\overbrace{\sum_{r=2}^{n}(-1)^{r-1}\coth\left(\frac{z\log r}{2}\right)}^{\text{$l=$ number of terms is even}}\right), & \text{if $n$ is even}\\
\frac{1}{1-2^{1-z}}\left(\frac{1}{2}+\frac{1}{2}\overbrace{\sum_{r=2}^{n}(-1)^{r-1}\coth\left(\frac{z\log r}{2}\right)}^{\text{$l=$ number of terms is odd}}\right), 
& \text{if $n$ is odd}.\\
\end{cases}
\end{eqnarray*}
\end{proof}
%%%%%%%%%%%%%%%%%%%%%%%%%%%%%%%%%%%%%%%%%%%%%%%%%%%%%%%%%%%%%%%%%%%%%%%%%%%%%%%%%%%%%%%%%%%%%%%%%%%%%%%%%%
By taking two terms from the above approximation 
\begin{eqnarray*}
1-\frac{1}{2^x-1}=0.
\end{eqnarray*}
Then, there are infinitely many complex solutions $ x=1+\frac{2i\pi n}{\log 2} $ where $n$ is an integer numbers. It is easy to obtain 
these roots by solving it directly.  $2^{x-1}=1$. From here, by using $\log z$ where $z$  is a complex number. One gets, $ x=1+\frac{2i\pi n}{\log 2} $ 
(Vinner!) \cite{18}. For sone of numerical calculatons of the limiting approximations for Riemann zeta function obtained above, see \cite{17}. 

%%%%%%%%%%%%%%%%%%%%%%%%%%
\section{Conclusion}\label{con}
To the best of author's knowledge this is a very new approach applied to the Riemann zeta function, which is not reported in the existing literature yet. 
The author thinks that theorems obtained above  are good approximations to the limit of the Riemann zeta function for sufficiently large $n$. The limit obtained here is a very promising for the non trivial zeros of  the Rieman zeta function but further study is needed to give mathematically rigorous proofs.
%%%%%%%%%%%%%%%%%%%%%%%%%%%%%%%%%%%%%%%%%%%%%%%%%%%%%%%%%%%%%%%%%%%%%%%%%%%%%%%%%%%%%%%%%%%%%%%%%%%%%%%%%%%%%%%%

\end{document}